\documentclass[12pt]{amsart}
\usepackage{amssymb}
\newtheorem{theorem}[subsection]{Theorem}

\newtheorem{corollary}[subsection]{Corollary}

\theoremstyle{remark}

\numberwithin{equation}{section}
\newcommand{\topten}{\circledcirc}
\newcommand{\Topten}{\raisebox{-1pt}{\Large$\circledcirc$}}
\newcommand{\ffrac}[2]{{\mbox{\large$\frac{#1}{#2}$}}}

\newcommand{\Hom}{{\mathrm{Hom}}}

\begin{document}
\title[The Joseph Ideal for the Classical Groups]
{The Uniqueness of the Joseph Ideal for the Classical Groups}
\author{Michael Eastwood}
\address{Department of Mathematics, University of Adelaide,
South Australia 5005}
\email{meastwoo@maths.adelaide.edu.au}
\author{Petr Somberg}
\address{Mathematical Institute, Charles University, Sokolovsk\'a 83, Praha,
\newline\indent Czech Republic}
\email{somberg@karlin.mff.cuni.cz}
\author{Vladim\'{\i}r Sou\v{c}ek}
\address{Mathematical Institute, Charles University, Sokolovsk\'a 83, Praha,
\newline\indent Czech Republic}
\email{soucek@karlin.mff.cuni.cz}
\thanks{This research was supported by the Australian Research Council, the 
grants GA\v{C}R 201/03/P137, GA\v{C}R 201/05/2117, and by the research grant 
LC505. The first author would also like to thank the Eduard \v{C}ech Center 
and the Mathematical Institute of Charles University for their excellent 
hospitality during the preparation of this article.}
\begin{abstract}
The Joseph ideal is a special ideal in the universal enveloping algebra of a
simple Lie algebra. For the classical groups, its uniqueness is equivalent to
the existence of tensors with special properties. The existence of these
tensors is usually concluded abstractly but here we present explicit formulae.
This allows a rather direct computation of a parameter used in defining the
Joseph ideal.
\end{abstract}
\renewcommand{\subjclassname}{\textup{2000} Mathematics Subject Classification}
\maketitle

\section{Introduction}
Let ${\mathfrak{g}}$ be a simple complex Lie algebra. Let
${\mathfrak{g}}\topten{\mathfrak{g}}$ denote the Cartan product of
${\mathfrak{g}}$ with itself, namely the unique irreducible component of 
${\mathfrak{g}}\otimes{\mathfrak{g}}$ whose highest weight is twice the 
highest weight of ${\mathfrak{g}}$ itself. Also, let us denote by 
$${\mathfrak{g}}\otimes{\mathfrak{g}}\in X\otimes Y\longmapsto X\topten Y\in
{\mathfrak{g}}\topten{\mathfrak{g}}\subset{\mathfrak{g}}\otimes{\mathfrak{g}}$$
the invariant projection onto the Cartan product. Let $[X,Y]$ denote the Lie
bracket and $\langle X,Y\rangle$ the Killing form on~${\mathfrak{g}}$. Then in
the full tensor algebra ${\bigotimes}{\mathfrak{g}}$ let us consider the
two-sided ideal $I_\lambda$ generated by elements of the form 
\begin{equation}\label{generators}\textstyle
X\otimes Y - X\topten Y - \ffrac12 [X,Y] - \lambda\langle X,Y\rangle\in
\bigotimes^2{\mathfrak{g}}\oplus{\mathfrak{g}}\oplus{\mathbb{C}},\quad
\forall\, X,Y\in{\mathfrak{g}}.\end{equation}
Let us denote by $A_\lambda$, the quotient algebra
$\bigotimes{\mathfrak{g}}/I_\lambda$.

\begin{theorem}\label{unique} For each complex simple Lie algebra not
isomorphic to\/~${\mathfrak{sl}}(2,{\mathbb{C}})$, there is precisely one value
of $\lambda$ for which $A_\lambda$ is infinite dimensional.
\end{theorem}
Usually, the special linear series of Lie algebras are excluded from this
theorem. This is owing to the geometric realisation of certain representations
usually used in defining the Joseph ideal~\cite{j}. However, as pointed out by
Braverman and Joseph~\cite{bj}, if the result is stated as above, then the
special linear algebras save for ${\mathfrak{sl}}(2,{\mathbb{C}})$ are
included. The main difficulty in proving Theorem~\ref{unique} is in showing
that, for all $\lambda$ save for a special value, the ideal is of finite
codimension. Braverman and Joseph~\cite{bj} present an abstract argument for
this in general, which they make explicit for the symplectic and special linear
algebras. They remark, however, that `in general such verification seems very
difficult'. This is what we accomplish for the orthogonal algebras. For
completeness and convenience, we also present the special tensors that may be
used in direct proofs of the symplectic and special linear cases. This covers
the classical groups. In fact, as pointed out by Gan and Savin~\cite{gs}, the
exceptional cases are somewhat easier. These authors also tackle the classical
groups by direct but still quite intricate means. 

The first author would like to thank Nolan Wallach for drawing his attention to
the Joseph ideal in response to a talk given at the University of California,
San Diego, on higher symmetries of the Laplacian~\cite{e}.

\section{The orthogonal case}\label{two}
We shall use index conventions for tensors as is standard in differential
geometry. More precisely, we use the abstract index notation of
Penrose~\cite{OT}. In particular $g^{ab}$ will denote the non-degenerate
quadratic form preserved by ${\mathfrak{so}}(n,{\mathbb{C}})$ and $g_{ab}$ its
inverse. We shall `raise and lower' indices without comment: so $X^a=g^{ab}X_b$ 
and $X_a=g_{ab}X^b$ where a repeated index denotes the invariant pairing 
between vectors and covectors.

\begin{theorem}\label{sotheorem}
For $\lambda\not=-\ffrac{n-4}{4(n-1)(n-2)}$ and $n\geq 5$, the two-sided ideal
in $\bigotimes{\mathfrak{so}}(n,{\mathbb{C}})$ generated by 
$$\textstyle X\otimes
Y - X\topten Y - \ffrac12 [X,Y] - \lambda\langle X,Y\rangle,\quad
\mbox{for }X,Y\in{\mathfrak{so}}(n,{\mathbb{C}})$$
contains ${\mathfrak{so}}(n,{\mathbb{C}})$, the first graded piece of 
$\bigotimes{\mathfrak{so}}(n,{\mathbb{C}})$.
\end{theorem}
\begin{proof}
If we identify ${\mathfrak{so}}(n,{\mathbb{C}})$ with skew tensors in the usual 
fashion, then the ideal is generated by tensors of the form 
$$V^{abcd}-(\topten V)^{abcd}-\ffrac12(V^a{}_b{}^{bd}-V^d{}_b{}^{ba})
+\lambda(n-2)V^{ab}{}_{ab}\;\mbox{ for }V^{abcd}=-V^{bacd}=-V^{abdc},$$
where $(\topten V)^{abcd}$ denotes the Cartan part of $V^{abcd}$. (There is, of
course, an explicit formula for $(\topten V)^{abcd}$ but we shall not need all
of it.) Consider the following tensor 
$$\begin{array}{rcl} S^{abcdef}&=&
2g^{af}g^{be}T^{cd}-2g^{ae}g^{bf}T^{cd}-2g^{cf}g^{de}T^{ab}+2g^{ce}g^{df}T^{ab}
\\ &&{}
+g^{ac}g^{be}T^{df}-g^{bc}g^{ae}T^{df}-g^{ad}g^{be}T^{cf}+g^{bd}g^{ae}T^{cf}\\
&&\enskip{}
-g^{ac}g^{bf}T^{de}+g^{bc}g^{af}T^{de}+g^{ad}g^{bf}T^{ce}-g^{bd}g^{af}T^{ce}\\
&&\quad{}
-g^{ac}g^{de}T^{bf}+g^{ad}g^{ce}T^{bf}+g^{bc}g^{de}T^{af}-g^{bd}g^{ce}T^{af}\\
&&\quad\enskip{}
+g^{ac}g^{df}T^{be}-g^{ad}g^{cf}T^{be}-g^{bc}g^{df}T^{ae}+g^{bd}g^{cf}T^{ae}.
\end{array}$$
for $T^{ab}=-T^{ba}$. It is immediate that $S^{abcdef}=-S^{cdabef}$ and 
readily verified that
$$Z^{abcdef}\equiv 
\ffrac13(S^{abcdef}+S^{abefcd})+
\ffrac16(S^{abcedf}-S^{abdecf}-S^{abcfde}+S^{abdfce})$$
is given by
$$\begin{array}{rcl}
Z^{abcdef}&=&2g^{ce}g^{df}T^{ab} - 2g^{de}g^{cf}T^{ab}\\
&&{}   -\frac12 g^{ac}g^{de}T^{bf} +\frac12 g^{ad}g^{ce}T^{bf}
     +\frac12 g^{bc}g^{de}T^{af} -\frac12 g^{bd}g^{ce}T^{af}\\ 
&&\enskip{}   +\frac12 g^{ac}g^{df}T^{be} -\frac12 g^{ad}g^{cf}T^{be}
     -\frac12 g^{bc}g^{df}T^{ae} +\frac12 g^{bd}g^{cf}T^{ae}\\
&&\quad{}   -\frac12 g^{ae}g^{cf}T^{bd} +\frac12 g^{ae}g^{df}T^{bc} 
     +\frac12 g^{be}g^{cf}T^{ad} -\frac12 g^{be}g^{df}T^{ac}\\
&&\quad\enskip{}   +\frac12 g^{af}g^{ce}T^{bd} -\frac12 g^{af}g^{de}T^{bc}
     -\frac12 g^{bf}g^{ce}T^{ad} +\frac12 g^{bf}g^{de}T^{ac}.
\end{array}$$
Generally, $S^{abcdef}\mapsto Z^{abcdef}$ is the formula for the 
${\mathfrak{sl}}(n,{\mathbb{C}})$-invariant projection
$$\begin{picture}(10,30)
\put(0,5){\line(1,0){10}}
\put(0,15){\line(1,0){10}}
\put(0,25){\line(1,0){10}}
\put(0,5){\line(0,1){20}}
\put(10,5){\line(0,1){20}}
\end{picture}
\raisebox{12pt}{$\;\otimes\;$}
\begin{picture}(10,30)
\put(0,5){\line(1,0){10}}
\put(0,15){\line(1,0){10}}
\put(0,25){\line(1,0){10}}
\put(0,5){\line(0,1){20}}
\put(10,5){\line(0,1){20}}
\end{picture}
\raisebox{12pt}{$\;\otimes\;$}
\begin{picture}(10,30)
\put(0,5){\line(1,0){10}}
\put(0,15){\line(1,0){10}}
\put(0,25){\line(1,0){10}}
\put(0,5){\line(0,1){20}}
\put(10,5){\line(0,1){20}}
\end{picture}
\raisebox{12pt}{$\;\longrightarrow\;$}
\begin{picture}(10,30)
\put(0,5){\line(1,0){10}}
\put(0,15){\line(1,0){10}}
\put(0,25){\line(1,0){10}}
\put(0,5){\line(0,1){20}}
\put(10,5){\line(0,1){20}}
\end{picture}
\raisebox{12pt}{$\;\otimes\;$}
\begin{picture}(25,30)
\put(0,5){\line(1,0){20}}
\put(0,15){\line(1,0){20}}
\put(0,25){\line(1,0){20}}
\put(0,5){\line(0,1){20}}
\put(10,5){\line(0,1){20}}
\put(20,5){\line(0,1){20}}
\end{picture}
\raisebox{12pt}{.}$$
In this case, however, the result is manifestly pure trace in $cdef$. It 
follows that under the further ${\mathfrak{so}}(n,{\mathbb{C}})$-invariant 
projection 
$$\begin{picture}(10,30)
\put(0,5){\line(1,0){10}}
\put(0,15){\line(1,0){10}}
\put(0,25){\line(1,0){10}}
\put(0,5){\line(0,1){20}}
\put(10,5){\line(0,1){20}}
\end{picture}
\raisebox{12pt}{$\;\otimes\;$}
\begin{picture}(10,30)
\put(0,5){\line(1,0){10}}
\put(0,15){\line(1,0){10}}
\put(0,25){\line(1,0){10}}
\put(0,5){\line(0,1){20}}
\put(10,5){\line(0,1){20}}
\end{picture}
\raisebox{12pt}{$\;\otimes\;$}
\begin{picture}(10,30)
\put(0,5){\line(1,0){10}}
\put(0,15){\line(1,0){10}}
\put(0,25){\line(1,0){10}}
\put(0,5){\line(0,1){20}}
\put(10,5){\line(0,1){20}}
\end{picture}
\raisebox{12pt}{$\;\longrightarrow\;$}
\begin{picture}(10,30)
\put(0,5){\line(1,0){10}}
\put(0,15){\line(1,0){10}}
\put(0,25){\line(1,0){10}}
\put(0,5){\line(0,1){20}}
\put(10,5){\line(0,1){20}}
\end{picture}
\raisebox{12pt}{$\;\otimes\;$}
\begin{picture}(25,30)
\put(0,5){\line(1,0){20}}
\put(0,15){\line(1,0){20}}
\put(0,25){\line(1,0){20}}
\put(0,5){\line(0,1){20}}
\put(10,5){\line(0,1){20}}
\put(20,5){\line(0,1){20}}
\end{picture}
\raisebox{12pt}{$\;\longrightarrow\;$}
\begin{picture}(10,30)
\put(0,5){\line(1,0){10}}
\put(0,15){\line(1,0){10}}
\put(0,25){\line(1,0){10}}
\put(0,5){\line(0,1){20}}
\put(10,5){\line(0,1){20}}
\end{picture}
\raisebox{12pt}{$\;\otimes\;$}
\begin{picture}(25,30)
\put(0,5){\line(1,0){20}}
\put(0,15){\line(1,0){20}}
\put(0,25){\line(1,0){20}}
\put(0,5){\line(0,1){20}}
\put(10,5){\line(0,1){20}}
\put(20,5){\line(0,1){20}}
\put(22,7){\makebox(0,0)[l]{$\circ$}}
\end{picture}$$
where $\circ$ denotes the trace-free part, the tensor $S^{abcdef}$ maps to
zero. But, for $n\geq 5$, this is the Cartan part in $cdef$. Its skew symmetry
$S^{abcdef}=-S^{cdabef}$ ensures that the Cartan part is also zero with respect
to $abcd$. Therefore, we may immediately reduce $S^{abcdef}$ in two different
ways with respect to the given ideal. We obtain, after a short calculation,
$$S^a{}_b{}^{bdef}=
(n-4)[g^{af}T^{de}-g^{ae}T^{df}+g^{de}T^{af}-g^{df}T^{ae}],$$
which is skew in $ad$. Therefore $S^{ab}{}_{ab}{}^{cd}=0$ and
$$S^{abcdef}\simeq \ffrac12(S^a{}_b{}^{bdef}-S^d{}_b{}^{baef})
=(n-4)[g^{af}T^{de}-g^{ae}T^{df}+g^{de}T^{af}-g^{df}T^{ae}].$$
Tracing over $de$ now gives 
$$g^{af}T^{de}-g^{ae}T^{df}+g^{de}T^{af}-g^{df}T^{ae}\simeq
(n-2)T^{af}.$$
Altogether, 
\begin{equation}\label{firstreduction}
S^{abcdef}\simeq(n-2)(n-4)T^{af}.\end{equation}
On the other hand, 
$$S^{abc}{}_d{}^{df}=-g^{af}T^{bc}+g^{bf}T^{ac}-2(n-2)g^{cf}T^{ab}
-(n-3)g^{ac}T^{bf}+(n-3)g^{bc}T^{af}$$
so
$$S^{abc}{}_d{}^{df}-S^{abf}{}_d{}^{dc}=
-(n-4)[g^{ac}T^{bf}-g^{bc}T^{af}-g^{af}T^{bc}+g^{bf}T^{ac}]$$
and
$$S^{abcd}{}_{cd}=2(n-1)(n-2)T^{ab}.$$
Therefore, 
$$S^{abcdef}\simeq-
\ffrac{n-4}2[g^{ac}T^{bf}-g^{bc}T^{af}-g^{af}T^{bc}+g^{bf}T^{ac}]-
2\lambda(n-1)(n-2)^2T^{af}.$$
But tracing over $bc$ gives
$$g^{ac}T^{bf}-g^{bc}T^{af}-g^{af}T^{bc}+g^{bf}T^{ac}\simeq
-(n-2)T^{af}$$
and so
\begin{equation}\label{secondreduction}
S^{abcdef}\simeq(n-2)[\ffrac{n-4}2-2\lambda(n-1)(n-2)]T^{af}.\end{equation}
Comparing (\ref{firstreduction}) with~(\ref{secondreduction}), we conclude that
$T^{ab}$ must be in the ideal unless we have 
$\lambda=-\ffrac{n-4}{4(n-1)(n-2)}$. This is exactly what we wanted to prove.
\end{proof}

\section{The symplectic case} 
\begin{theorem}\label{sptheorem}
For $\lambda\not=-\ffrac{1}{16(n+1)}$ and~$n\geq 2$, the two-sided ideal
in $\bigotimes{\mathfrak{sp}}(2n,{\mathbb{C}})$ generated by 
$$\textstyle X\otimes
Y - X\topten Y - \ffrac12 [X,Y] - \lambda\langle X,Y\rangle,\quad
\mbox{for }X,Y\in{\mathfrak{sp}}(2n,{\mathbb{C}})$$
contains ${\mathfrak{sp}}(2n,{\mathbb{C}})$, the first graded piece of 
$\bigotimes{\mathfrak{sp}}(2n,{\mathbb{C}})$.
\end{theorem}
\begin{proof}
Let $\omega^{ab}$ denote the skew form preserved by
${\mathfrak{sp}}(2n,{\mathbb{C}})$ and adopt the convention that
$\omega^{ac}\omega_{bc}$ is the identity. In particular
$\omega^{ab}\omega_{ab}=2n$. If we use $\omega_{ab}$ to lower indices according
to $X_b=X^a\omega_{ab}$, then we may identify
${\mathfrak{sp}}(2n,{\mathbb{C}})$ as symmetric tensors $T^{ab}=T^{ab}$ and the
ideal is generated by
$$V^{abcd}-(\topten V)^{abcd}-\ffrac12(V^a{}_b{}^{bd}+V^d{}_b{}^{ba})
+2\lambda(n+1)V^{ab}{}_{ab}\;\mbox{ for }V^{abcd}=V^{bacd}=V^{abdc}.$$
Now consider the tensor 
$$\begin{array}{rcl}
S^{abcdef}&=& 4\omega^{af}\omega^{be}T^{cd}+4\omega^{ae}\omega^{bf}T^{cd}
-4\omega^{cf}\omega^{de}T^{ab}-4\omega^{ce}\omega^{df}T^{ab} \\ &&{}
-\omega^{ac}\omega^{be}T^{df}-\omega^{bc}\omega^{ae}T^{df}
-\omega^{ad}\omega^{be}T^{cf}-\omega^{bd}\omega^{ae}T^{cf}\\
&&\enskip{}
-\omega^{ac}\omega^{bf}T^{de}-\omega^{bc}\omega^{af}T^{de}
-\omega^{ad}\omega^{bf}T^{ce}-\omega^{bd}\omega^{af}T^{ce}\\
&&\quad{}
-\omega^{ac}\omega^{de}T^{bf}-\omega^{ad}\omega^{ce}T^{bf}
-\omega^{bc}\omega^{de}T^{af}-\omega^{bd}\omega^{ce}T^{af}\\
&&\quad\enskip{}
-\omega^{ac}\omega^{df}T^{be}-\omega^{ad}\omega^{cf}T^{be}
-\omega^{bc}\omega^{df}T^{ae}-\omega^{bd}\omega^{cf}T^{ae}.
\end{array}$$
for $T^{ab}=T^{ba}$. It is immediate that $S^{abcdef}=-S^{cdabef}$ and 
readily verified that     
$$Z^{abcdef}\equiv\ffrac16(
S^{abcdef}+S^{abdefd}+S^{abcfde}+S^{abfedc}+S^{abedfc}+S^{abdfec})$$
vanishes. This is already the Cartan part with respect to the $cdef$ indices.
Therefore, we may reduce $S^{abdcef}$ modulo the ideal in two
different ways. We obtain
$$S^{abcdef}\simeq -2(n-1)(\omega^{ae}T^{df}+\omega^{de}T^{af}
+\omega^{af}T^{de}+\omega^{df}T^{ae})\simeq-4(n-1)(n+1)T^{af}$$
or
$$\begin{array}lS^{abcdef}\simeq
-(n-1)(\omega^{ac}T^{bf}+\omega^{bc}T^{af}
+\omega^{af}T^{bc}+\omega^{bf}T^{ac})+32\lambda(n-1)(n+1)^2T^{af}\\[3pt]
\phantom{S^{abcdef}}\simeq-2(n-1)(n+1)T^{af}+32\lambda(n-1)(n+1)^2T^{af}.
\end{array}$$
Comparing these two reductions, we see that $T^{ab}$ lies in the ideal unless
$\lambda=-\ffrac1{16(n+1)}$. This is what we wanted to prove. 
\end{proof}
Note that the critical value of $\lambda$ for 
${\mathfrak{sp}}(4,{\mathbb{C}})\cong{\mathfrak{so}}(5,{\mathbb{C}})$
may be computed either from Theorem~\ref{sotheorem} or 
Theorem~\ref{sptheorem}. Its common value is $-1/48$. 

\section{The special linear case}\label{four}
\begin{theorem}\label{sltheorem}
For $\lambda\not=-\ffrac{1}{8(n+1)}$ and~$n\geq 3$, the two-sided ideal
in $\bigotimes{\mathfrak{sl}}(n,{\mathbb{C}})$ generated by 
$$\textstyle X\otimes
Y - X\topten Y - \ffrac12 [X,Y] - \lambda\langle X,Y\rangle,\quad
\mbox{for }X,Y\in{\mathfrak{sl}}(n,{\mathbb{C}})$$
contains ${\mathfrak{sl}}(n,{\mathbb{C}})$, the first graded piece of 
$\bigotimes{\mathfrak{sl}}(n,{\mathbb{C}})$.
\end{theorem}
\begin{proof}
If we identify ${\mathfrak{sl}}(n,{\mathbb{C}})$ with trace-free tensors 
$T^a{}_b$ in the usual manner, then the ideal is generated by tensors of the 
form 
\begin{equation}\label{slideal}
V^a{}_b{}^c{}_d-(\topten V)^a{}_b{}^c{}_d
-\ffrac12(V^a{}_b{}^b{}_d-V^b{}_d{}^a{}_b)-2\lambda nV^a{}_b{}^b{}_a
\quad\mbox{for }V^a{}_a{}^c{}_d=0=V^a{}_b{}^c{}_c.
\end{equation}
Consider the tensor    
$$\begin{array}{rcl}
S^a{}_b{}^c{}_d{}^e{}_f&=& 
 \delta^e{}_d\delta^c{}_fT^a{}_b-\frac1{n}\delta^c{}_d\delta^e{}_fT^a{}_b
-\delta^e{}_b\delta^a{}_fT^c{}_d
+\frac1{n}\delta^a{}_b\delta^e{}_fT^c{}_d\\[3pt] 
&&{}
+\delta^a{}_d\delta^e{}_bT^c{}_f-\frac1{n}\delta^a{}_d\delta^e{}_fT^c{}_b
-\delta^c{}_b\delta^e{}_dT^a{}_f+\frac1{n}\delta^c{}_b\delta^e{}_fT^a{}_d
\end{array}$$
for $T^a{}_a=0$. It is immediate that
$S^a{}_b{}^c{}_d{}^e{}_f=-S^c{}_d{}^a{}_b{}^e{}_f$. In particular, the
Cartan part of $S^a{}_b{}^c{}_d{}^e{}_f$ with respect to the indices $abcd$
vanishes. Hence, we may use (\ref{slideal}) to reduce
$S^a{}_b{}^c{}_d{}^e{}_f$ modulo the given ideal. We obtain
$$S^a{}_b{}^c{}_d{}^e{}_f\simeq 
-(n-1)\delta^e{}_dT^a{}_f-\delta^a{}_fT^e{}_d+\delta^a{}_dT^e{}_f
+\delta^e{}_fT^a{}_d\simeq-\ffrac12n(n-2)T^a{}_f.$$
On the other hand, it is readily verified that 
$$Z^a{}_b{}^c{}_d{}^e{}_f\equiv\ffrac14(
S^a{}_b{}^c{}_d{}^e{}_f+S^a{}_b{}^e{}_d{}^c{}_f+S^a{}_b{}^c{}_f{}^e{}_d
+S^a{}_b{}^e{}_f{}^c{}_d)$$
is given by 
$$\begin{array}{rcl}Z^a{}_b{}^c{}_d{}^e{}_f&=&
\frac12\delta^e{}_d\delta^c{}_fT^a{}_b
-\frac14\delta^c{}_b\delta^e{}_dT^a{}_f
-\frac1{2n}\delta^c{}_d\delta^e{}_fT^a{}_b
+\frac1{4n}\delta^a{}_b\delta^e{}_fT^c{}_d\\[3pt]
&&{}-\frac1{4n}\delta^a{}_d\delta^e{}_fT^c{}_b
+\frac1{4n}\delta^c{}_b\delta^e{}_fT^a{}_d
-\frac14\delta^e{}_b\delta^c{}_fT^a{}_d
{}-\frac1{4n}\delta^c{}_d\delta^a{}_fT^e{}_b\\[3pt]
&&\enskip{}+\frac1{4n}\delta^a{}_b\delta^e{}_dT^c{}_f
+\frac1{4n}\delta^a{}_b\delta^c{}_fT^e{}_d
+\frac1{4n}\delta^c{}_d\delta^e{}_bT^a{}_f
-\frac1{2n}\delta^c{}_f\delta^e{}_dT^a{}_b\\[3pt]
&&\quad{}+\frac1{4n}\delta^a{}_b\delta^c{}_dT^e{}_f
-\frac1{4n}\delta^a{}_f\delta^e{}_dT^c{}_b
+\frac1{4n}\delta^c{}_b\delta^e{}_dT^a{}_f
-\frac1{4n}\delta^a{}_d\delta^c{}_fT^e{}_b\\[3pt]
&&\enskip\quad{}+\frac1{4n}\delta^e{}_b\delta^c{}_fT^a{}_d
+\frac12\delta^e{}_f\delta^c{}_dT^a{}_b
-\frac14\delta^c{}_b\delta^e{}_fT^a{}_d
-\frac14\delta^e{}_b\delta^c{}_dT^a{}_f.
\end{array}$$
Generally, $S^a{}_b{}^c{}_d{}^e{}_f\mapsto Z^a{}_b{}^c{}_d{}^e{}_f$ followed by
the removal of all traces in the $cdef$ indices is the Cartan projection in
these indices. In this case, however, $Z^a{}_b{}^c{}_d{}^e{}_f$ is manifestly
pure trace and so this Cartan part of $S^a{}_b{}^c{}_d{}^e{}_f$ vanishes. This 
allows us to use (\ref{slideal}) with respect to the $cdef$ indices to conclude
that 
$$\begin{array}{rcl}S^a{}_b{}^c{}_d{}^e{}_f&\simeq&\ffrac12(\delta^c{}_fT^a{}_b
-(n-1)\delta^c{}_bT^a{}_f+\delta^a{}_bT^c{}_f-\delta^a{}_fT^c{}_b)
+2\lambda n(n-2)(n+1)T^a{}_f\\[3pt]
&\simeq&-\ffrac14n(n-2)T^a{}_f+2\lambda n(n-2)(n+1)T^a{}_f.
\end{array}$$
Comparing this with our previous reduction, we see that $T^a{}_b$ lies in the 
ideal unless $\lambda=-\ffrac1{8(n+1)}$. This is what we wanted to prove.
\end{proof}
Note that the critical value of $\lambda$ for
${\mathfrak{sl}}(4,{\mathbb{C}})\cong{\mathfrak{so}}(6,{\mathbb{C}})$ may be
computed either from Theorem~\ref{sotheorem} or Theorem~\ref{sltheorem}. Its
common value is $-1/40$. 

\section{Remarks and Conclusions}
We should explain how the special tensors used in the proofs of
Theorems~\ref{sotheorem}, \ref{sptheorem}, and~\ref{sltheorem} arise. For all
simple Lie algebras other than the special linear series, there is a common
source as follows. Let ${\mathfrak{g}}$ denote a complex simple Lie algebra
and let $\Phi$ denote the composition
$$\Lambda^2{\mathfrak{g}}\otimes{\mathfrak{g}}\hookrightarrow
{\mathfrak{g}}\otimes{\mathfrak{g}}\otimes{\mathfrak{g}}
\xrightarrow{{\mathrm{Id}}\,\otimes\underbar{\phantom{x}}\topten
\underbar{\phantom{x}}}
{\mathfrak{g}}\otimes\Topten^2{\mathfrak{g}}.$$
\begin{theorem} For any simple complex Lie algebra~${\mathfrak{g}}$ not
isomorphic to ${\mathfrak{sl}}(n,{\mathbb{C}})$, 
$$\dim\Hom_{\mathfrak{g}}({\mathfrak{g}},
\Lambda^2{\mathfrak{g}}\otimes{\mathfrak{g}})=2\quad\mbox{and}\quad
\dim\Hom_{\mathfrak{g}}({\mathfrak{g}},
{\mathfrak{g}}\otimes\Topten^2{\mathfrak{g}})=1.$$
\end{theorem}
\begin{proof} A case-by-case verification using, for example, Klimyk's formula.
\end{proof}
\begin{corollary} For any simple complex Lie algebra~${\mathfrak{g}}$ not
isomorphic to ${\mathfrak{sl}}(n,{\mathbb{C}})$,
$$\dim\Hom_{\mathfrak{g}}({\mathfrak{g}},\ker\Phi)\geq 1.$$
\end{corollary}
\noindent This result is used abstractly by Braverman and Joseph~\cite{bj} and
our proofs are very much motivated by this approach: in proving
Theorems~\ref{sotheorem} and~\ref{sptheorem} we find explicit non-zero
homomorphisms ${\mathfrak{g}}\to\ker\Phi$. In fact, it is easily verified that
$\Hom_{\mathfrak{g}}({\mathfrak{g}},\ker\Phi)$ is 1-dimensional so our
homomorphisms are unique up to scale. 

For the special linear algebras the dimensions are different:--
\begin{theorem} For ${\mathfrak{g}}={\mathfrak{sl}}(n,{\mathbb{C}})$, 
$$\dim\Hom_{\mathfrak{g}}({\mathfrak{g}},
\Lambda^2{\mathfrak{g}}\otimes{\mathfrak{g}})=
\left\{\begin{array}{rl}4&\mbox{if }n\geq3\\ 1&\mbox{if }n=2
\end{array}\right.
\mbox{ and}\quad
\dim\Hom_{\mathfrak{g}}({\mathfrak{g}},
{\mathfrak{g}}\otimes\Topten^2{\mathfrak{g}})=1,\;\forall n\geq 2.$$
\end{theorem}
\begin{corollary}  For ${\mathfrak{g}}={\mathfrak{sl}}(n,{\mathbb{C}})$ with
$n\geq 3$, 
$$\dim\Hom_{\mathfrak{g}}({\mathfrak{g}},\ker\Phi)\geq 3.$$
\end{corollary}
\noindent In fact, using tensors, we have checked that 
$\Hom_{\mathfrak{g}}({\mathfrak{g}},\ker\Phi)$ is 3-dimensional and within it
there is a 2-dimensional subspace
$\Hom_{\mathfrak{g}}({\mathfrak{g}},\ker\Phi\cap\ker\Psi)$ where $\Psi$ is the 
composition
$$\Lambda^2{\mathfrak{g}}\otimes{\mathfrak{g}}\hookrightarrow
{\mathfrak{g}}\otimes{\mathfrak{g}}\otimes{\mathfrak{g}}
\xrightarrow{{\mathrm{Id}}\,\otimes\langle\underbar{\phantom{x}},
\underbar{\phantom{x}}\rangle}
{\mathfrak{g}}\otimes{\mathbb{C}}={\mathfrak{g}}.$$
Any homomorphism in $\Hom_{\mathfrak{g}}({\mathfrak{g}},\ker\Phi)\setminus
\Hom_{\mathfrak{g}}({\mathfrak{g}},\ker\Phi\cap\ker\Psi)$ will suffice for
deriving the critical value of $\lambda$ as in our proof of
Theorem~\ref{sltheorem}. This critical value of $\lambda$ is also obtained by
Braverman and Joseph~\cite[\S7.4~and~\S7.7]{bj}. They also remark
\cite[\S5.4]{bj} that the symplectic case may be dealt with by an `extremely
rare' but `simple-minded procedure' going back to Dirac. {From} the tensorial
point of view, the reason for this is that if one na\"{\i}vely extends a tensor
$$S\in\big(\Lambda^2{\mathfrak{sp}}(2m,{\mathbb{C}})\otimes
{\mathfrak{sp}}(2m,{\mathbb{C}})\big)\cap 
\big({\mathfrak{sp}}(2m,{\mathbb{C}})\otimes
\Topten^2{\mathfrak{sp}}(2m,{\mathbb{C}})\big)$$
by adding zero components then one obtains a tensor in 
$$\big(\Lambda^2{\mathfrak{sp}}(2n,{\mathbb{C}})\otimes
{\mathfrak{sp}}(2n,{\mathbb{C}})\big)\cap 
\big({\mathfrak{sp}}(2n,{\mathbb{C}})\otimes
\Topten^2{\mathfrak{sp}}(2n,{\mathbb{C}})\big)$$
for any $n>m$. In effect, Braverman and Joseph use this observation and an
explicit tensor for the case $n=2$ to obtain the general case. 

Usually, Theorem~\ref{unique} is stated in terms of the universal enveloping
algebra ${\mathfrak{U}}({\mathfrak{g}})$ of~${\mathfrak{g}}$. To do this,
notice that the generators (\ref{generators}) of $I_\lambda$ may be split into 
skew and symmetric parts:--
$$X\otimes Y -Y\otimes X - [X,Y]\quad\mbox{and}\quad
  X\otimes Y +Y\otimes X - 2X\topten Y - 2\lambda\langle X,Y\rangle$$
and that we may define the algebra $A_\lambda$ in two steps, firstly taking 
the quotient of the tensor algebra by the skew generators. This gives
${\mathfrak{U}}({\mathfrak{g}})$ and an image ideal $\bar I_\lambda$ so that
$A_\lambda={\mathfrak{U}}({\mathfrak{g}})/\bar I_\lambda$. What we have shown
more precisely in \S\ref{two}--\S\ref{four} is the following:--
\begin{theorem} For the classical complex simple algebras
$$\begin{array}{rcl}
{\mathfrak{g}}={\mathfrak{so}}(n,{\mathbb{C}}),&n\geq 5,&
\lambda\not =-\ffrac{n-4}{4(n-1)(n-2)}\\
{\mathfrak{g}}={\mathfrak{sp}}(2n,{\mathbb{C}}),&n\geq 2,&
\lambda\not=-\ffrac{1}{16(n+1)}\\
{\mathfrak{g}}={\mathfrak{sl}}(n,{\mathbb{C}}),&n\geq 3,&
\lambda\not=-\ffrac{1}{8(n+1)}
\end{array}$$
the ideal $\bar I_\lambda$ coincides with ${\mathfrak{U}}({\mathfrak{g}})$ if
$\lambda\not=0$ whilst $\bar I_0=
{\mathfrak{U}}_+({\mathfrak{g}})\subset{\mathfrak{U}}({\mathfrak{g}})$, the
unique maximal ideal consisting of elements without constant part. 
\end{theorem}
\begin{proof}
Theorems~\ref{sotheorem}, \ref{sptheorem}, and~\ref{sltheorem} say that, in
these circumstances, the ideal $I_\lambda$ contains ${\mathfrak{g}}$ and hence
contains $\bigoplus_{s\geq 1}\bigotimes^s\!{\mathfrak{g}}$, whose image is
${\mathfrak{U}}_+({\mathfrak{g}})$ by definition. The conclusions are now
immediate from~(\ref{generators}).
\end{proof}
In all other cases the algebra $A_\lambda$ is, in fact, infinite-dimensional.
For the orthogonal algebras, for example, there are linear differential
operators 
$${\mathcal{D}}_X\quad\mbox{for all }X\in\Topten^s{\mathfrak{so}}(m+1,1)$$ 
constructed in \cite{e} that satisfy
$${\mathcal{D}}_X{\mathcal{D}}_Y={\mathcal{D}}_{X\topten Y}
+\ffrac12{\mathcal{D}}_{[X,Y]}-\ffrac{m-2}{4m(m+1)}{\mathcal{D}}_{\langle
X,Y\rangle},\quad\forall X,Y\in{\mathfrak{so}}(m+1,1).$$
The corresponding holomorphic differential operators provide a realisation of
$A_\lambda$ for ${\mathfrak{g}}={\mathfrak{so}}(n,{\mathbb{C}})$ and
$\lambda=-\ffrac{(n-4)}{4(n-1)(n-2)}$. There are similar linear
holomorphic differential operators for
${\mathfrak{g}}={\mathfrak{sl}}(n,{\mathbb{C}})$ constructed as follows. Recall
that in~\S\ref{four} we identified ${\mathfrak{sl}}(n,{\mathbb{C}})$ with
trace-free tensors~$X^a{}_b$. More generally,
$$\Topten^s{\mathfrak{sl}}(n,{\mathbb{C}})=\left\{
X\underbrace{{}^a{}_b{}^c{}_d{}^{\cdots}{}_{\cdots}{}^e{}_f}_{
2s\mbox{ \scriptsize indices}}
\quad\mbox{s.t. }
\begin{array}l
X{}^a{}_b{}^c{}_d{}^{\cdots}{}_{\cdots}{}^e{}_f=
X{}^{(a}{}_{(b}{}^c{}_d{}^{\cdots}{}_{\cdots}{}^{e)}{}_{f)}\\[3pt]
X{}^a{}_b{}^c{}_d{}^{\cdots}{}_{\cdots}{}^e{}_f\mbox{ is totally trace-free}.
\end{array}\right\}$$
and we define
$${\mathcal{D}}_X\equiv(-1)^sX{}^a{}_b{}^c{}_d{}^{\cdots}{}_{\cdots}{}^e{}_f
Z^bZ^d\cdots Z^f
\frac{\partial^s}{\partial Z^a\partial Z^c\cdots\partial Z^e}$$
as a holomorphic differential operator acting on~${\mathbb{C}}^n$. For
$X^a{}_b,Y^c{}_d\in{\mathfrak{sl}}(n,{\mathbb{C}})$,
$${\mathcal{D}}_X{\mathcal{D}}_Y-{\mathcal{D}}_Y{\mathcal{D}}_X
=(Y^a{}_cX^c{}_b-X^a{}_cY^c{}_b)Z^b\frac{\partial}{\partial Z^a}
=-[X,Y]^a{}_bZ^b\frac{\partial}{\partial Z^a}={\mathcal{D}}_{[X,Y]}$$
and 
$${\mathcal{D}}_X{\mathcal{D}}_Y+{\mathcal{D}}_Y{\mathcal{D}}_X
=2X^{(a}{}_{(b}Y^{c)}{}_{d)}Z^bZ^d\frac{\partial^2}{\partial Z^a\partial Z^c}
+(X^a{}_cY^c{}_b+Y^a{}_cX^c{}_b)Z^b\frac{\partial}{\partial Z^a}.$$
However, if we write 
$$X^{(a}{}_{(b}Y^{c)}{}_{d)}=C^{ac}{}_{bd}+D^{(a}{}_{(b}\delta^{c)}{}_{d)}
+E\delta^{(a}{}_{(b}\delta^{c)}{}_{d)},$$
where 
$$D^a{}_b=\ffrac1{n+2}(X^a{}_cY^c{}_b+Y^a{}_cX^c{}_b)
-\ffrac2{n(n+2)}X^c{}_dY^d{}_c\delta^a{}_b\quad\mbox{and}\quad
E=\ffrac1{n(n+1)}X^c{}_dY^d{}_c,$$
then $C^{ac}{}_{bd}$ and $D^a{}_b$ are trace-free. In particular,
$C^{ac}{}_{bd}=(X\topten Y)^{ac}{}_{bd}$ and
$$\begin{array}{rcl}
{\mathcal{D}}_X{\mathcal{D}}_Y+{\mathcal{D}}_Y{\mathcal{D}}_X
&=&\displaystyle 2(C^{ac}{}_{bd}+D^c{}_b\delta^a{}_d
+E\delta^c{}_b\delta^a{}_d)
Z^bZ^d\frac{\partial^2}{\partial Z^a\partial Z^c}\\
&&\displaystyle\quad{}
+(X^a{}_cY^c{}_b+Y^a{}_cX^c{}_b)Z^b\frac{\partial}{\partial Z^a}\\
&=&\displaystyle 2{\mathcal{D}}_{X\topten Y}
+2D^a{}_bZ^bZ^c\frac{\partial^2}{\partial Z^c\partial Z^a}
+2EZ^bZ^c\frac{\partial^2}{\partial Z^c\partial Z^b}\\
&&\displaystyle\quad{}
+(X^a{}_cY^c{}_b+Y^a{}_cX^c{}_b)Z^b\frac{\partial}{\partial Z^a}.
\end{array}$$
Now, let us restrict the action of these differential operators to germs $\phi$
of holomorphic functions defined near some basepoint in
${\mathbb{C}}^n\setminus\{0\}$ and `homogeneous of degree~$w$' in the sense
that $Z^a\partial/\partial Z^a\phi=w\phi$. We find that
$$\begin{array}{rcl}
{\mathcal{D}}_X{\mathcal{D}}_Y+{\mathcal{D}}_Y{\mathcal{D}}_X 
&=&\displaystyle
2{\mathcal{D}}_{X\topten Y}
+(2\ffrac{w-1}{n+2}+1)(X^a{}_cY^c{}_b+Y^a{}_cX^c{}_b)
Z^b\frac{\partial}{\partial Z^a}\\
&&\displaystyle\quad{}
+2(\ffrac1{n(n+1)}-\ffrac2{n(n+2)})w(w-1)X^c{}_dY^d{}_c\\
&=&\displaystyle
2{\mathcal{D}}_{X\topten Y}
+\ffrac{2w+n}{n+2}(X^a{}_cY^c{}_b+Y^a{}_cX^c{}_b)
Z^b\frac{\partial}{\partial Z^a}\\
&&\displaystyle\quad{}
-2\ffrac{w(w-1)}{(n+1)(n+2)}X^c{}_dY^d{}_c.
\end{array}$$
Assembling these computations we conclude that 
$${\mathcal{D}}_X{\mathcal{D}}_Y=
{\mathcal{D}}_{X\topten Y}+\ffrac{2w+n}{2(n+2)}(X^a{}_cY^c{}_b+Y^a{}_cX^c{}_b)
Z^b\frac{\partial}{\partial Z^a}+\ffrac12{\mathcal{D}}_{[X,Y]}-
\ffrac{w(w-1)}{2n(n+1)(n+2)}{\mathcal{D}}_{\langle X,Y\rangle}.$$
In particular, for $w=-n/2$ we obtain
$${\mathcal{D}}_X{\mathcal{D}}_Y=
{\mathcal{D}}_{X\topten Y}+\ffrac12{\mathcal{D}}_{[X,Y]}-
\ffrac1{8(n+1)}{\mathcal{D}}_{\langle X,Y\rangle}.$$
These operators provide a realisation of A$_\lambda$ for
${\mathfrak{g}}={\mathfrak{sl}}(n,{\mathbb{C}})$ and $\lambda=-\ffrac1{8(n+1)}$
and, in particular, show that this algebra is infinite dimensional.

If $n=2$ we can proceed further because, in this case,
$$X^a{}_cY^c{}_b+Y^a{}_cX^c{}_b=X^d{}_cY^c{}_d\delta^a{}_b
=\ffrac14\langle X,Y\rangle\delta^a{}_b$$
whence
$${\mathcal{D}}_X{\mathcal{D}}_Y=
{\mathcal{D}}_{X\topten Y}+\ffrac12{\mathcal{D}}_{[X,Y]}
+\ffrac{w(w+2)}{24}{\mathcal{D}}_{\langle X,Y\rangle}$$
for any $w\in{\mathbb{C}}$. In particular, this shows that $A_\lambda$ is
infinite-dimensional for ${\mathfrak{sl}}(2,{\mathbb{C}})$ no matter what
is~$\lambda$. 

An alternative to these geometric realisations of $A_\lambda$ is provided by
the generalised Poincar\'e-Birkhoff-Witt Theorem of Braverman and
Gaitsgory~\cite{bg}, which enables one to identify the associated graded
algebra~${\mathrm{gr}}(A_\lambda)$. Specifically, if we let
$R\subset{\mathfrak{g}}\otimes{\mathfrak{g}}$ be the ${\mathfrak{g}}$-invariant
complement to ${\mathfrak{g}}\topten{\mathfrak{g}}$ and $J(R)$ be the two-sided
ideal in $\bigotimes{\mathfrak{g}}$ generated by~$R$, then as a special case of
\cite{bg} we obtain criteria under which the canonical surjection
$p:\bigotimes{\mathfrak{g}}/J(R)\to{\mathrm{gr}}(A_\lambda)$ of graded algebras
is an isomorphism. These criteria are then verified by Braverman and
Joseph~\cite{bj} in the case of critical~$\lambda$. It follows from a result of
Kostant (given in a lecture at MIT in 1980 and explained with proof
in~\cite[Chapter~3]{g}) that the graded algebra $\bigotimes{\mathfrak{g}}/J(R)$
is simply the Cartan algebra
$\Topten{\mathfrak{g}}=\bigoplus_{s=0}^\infty\Topten^s{\mathfrak{g}}$ for any
complex simple Lie algebra~${\mathfrak{g}}$. This is also proved by tensorial
means in~\cite{e} for the orthogonal algebras, in~\cite{iciam} for the special
linear algebras, and the symplectic algebras are easily dealt with by a similar
argument. In~\cite{iciam}, however, it was incorrectly asserted that $p$ is
always an isomorphism.

\end{document}